 \documentclass[12pt]{article}
\usepackage{amssymb}
\newtheorem{theorem}{Theorem}
\newenvironment{proof}{\begin{trivlist}\item[]
{\bf Proof.} }{\hfill $\square$ \end{trivlist}}
\begin{document}
\begin{center}
{\LARGE\bf Continuous Invariants of Isolated}\\[10pt]
{\LARGE\bf Hypersurface Singularities}\\[15pt]
\end{center}
\begin{flushleft}
{\large\bf Michael G. Eastwood}\\
{\small Pure Mathematics Department, Adelaide University,
        South Australia 5005}\\
{\small E-mail: meastwoo@maths.adelaide.edu.au}
\end{flushleft}
\def\thefootnote{}
\footnotetext{This research was supported by the Australian Research Council.}
\def\thefootnote{$\sharp$}
\section{Introduction}
Suppose $V$ is a complex hypersurface in ${\mathbb C}^n$ with an isolated
singularity at the origin.  Writing $V$ as the zero locus of a holomorphic
function $f(z_1,\dots,z_n)$, it is well-known that the `moduli algebra'
$$A={\mathbb C}\{z_1,\dots,z_n\}/
(f,\partial f/\partial z_1,\dots,\partial f/\partial z_n)$$
is a finite-dimensional algebra depending only on the germ of $V$ at the origin
and is invariant under holomorphic change of co\"ordinates.  Remarkably, it was
shown by Mather and Yau~\cite{my} that $A$ completely determines the germ
of~$V$.  Thus, one should be able to distinguish between biholomorphically
inequivalent singularities on the basis of their corresponding moduli algebras.
This was accomplished for the simple elliptic singularities $\tilde E_7$ and
$\tilde E_8$ by Seeley and Yau \cite{sy} and for $\tilde E_6$ by Chen, Seeley,
and Yau~\cite{csy}.  In particular, they recovered Saito's
$j$-invariant~\cite{s} for $\tilde E_7$ and $\tilde E_8$ and detected an error
in his calculation for $\tilde E_6$.  In~\cite{s}, the $j$-invariant was
defined as that associated to the elliptic curve
arising as the exceptional divisor of the minimal resolution of~$V$.

The analysis in \cite{csy,sy} is fairly involved and specific to $\tilde E_6$,
$\tilde E_7$, and~$\tilde E_8$.  In this article we show how to employ
classical invariant theory to obtain $j$ directly and rather quickly from~$A$.
We explain $\tilde E_6$, $\tilde E_7$, and $\tilde E_8$ in detail but the
approach is quite general.  To illustrate this we present another example.

I would like to thank Stephen Yau for drawing my attention to this problem and
for several interesting conversations.  Thanks are also due to Adam Harris for
further useful conversations.

\section{Constructing Invariants}
The moduli algebra $A$ evidently has a unique maximal ideal ${\mathfrak m}$
generated by the co\"ordinate functions $z_1,\dots,z_n$.  Suppose that
${\mathfrak m}^N$ is one-dimensional and ${\mathfrak m}^{N+1}=0$.  Then the
algebra multiplication defines a linear transformation
$\bigodot^N{\mathfrak m}/{\mathfrak m}^2\to{\mathfrak m}^N$ or, in other words,
a tensor $a\in\bigodot^N({\mathfrak m}/{\mathfrak m}^2)^*$, canonically defined
up to scale.  The special linear group
${\mathrm {SL}}({\mathfrak m}/{\mathfrak m}^2)$ acts on $a$ and we may consider
its classical invariants.  By definition these are the polynomials in the
coefficients of $a$ that are invariant under the action of
${\mathrm {SL}}({\mathfrak m}/{\mathfrak m}^2)$.  Weyl's classical invariant
theory \cite{w} dictates the general form of such invariants.  Some discussion
and examples are provided in an appendix.
\begin{theorem}\label{construction} Suppose $P$ and $Q$ are classical
invariants of $a$ and homogeneous of the same degree.  Suppose $Q$ is non-zero.
Then $P/Q$ is an absolute invariant of the algebra~$A$.  More precisely, if $B$
is complex algebra isomorphic to $A$ then it must also be local, say with
maximal ideal~${\mathfrak n}$ of the same dimension as ${\mathfrak m}$, also
with ${\mathfrak n}^N$ one-dimensional and ${\mathfrak n}^{N+1}=0$ and so that
the value of $P/Q$ computed for the multiplication tensor $\bigodot^N{\mathfrak
n}/{\mathfrak n}^2\to{\mathfrak n}^N$ is the same as that computed for~$A$.
\end{theorem}

\begin{proof} Let us write $\rho$ for the representation of
${\mathrm {GL}}({\mathfrak m}/{\mathfrak m}^2)$ on
$\bigodot^N({\mathfrak m}/{\mathfrak m}^2)^*$.  To say that $P(a)$ is a
classical invariant of a tensor
$a\in\bigodot^N({\mathfrak m}/{\mathfrak m}^2)^*$ is to say that
$P(a)=P(\rho(M)a)$ for any $M\in{\mathrm {SL}}({\mathfrak m}/{\mathfrak m}^2)$.
Let $n=\dim_{\mathbb C}({\mathfrak m}/{\mathfrak m}^2)$.  If $P$ is homogeneous
of degree~$d$, then according to the classical theory (explained in the
appendix), $n$ must divide $dN$ and $P(a)=\det(M)^{dN/n}P(\rho(M)a)$
for $M\in{\mathrm {GL}}({\mathfrak m}/{\mathfrak m}^2)$.  Therefore, the ratio
$P(a)/Q(a)$ does not see any choice of basis.  Neither does it see the scale
ambiguity in~$a$: replacing $a$ by $\lambda a$ scales both $P$ and $Q$
by~$\lambda^d$.  These were the only arbitrary ingredients in the computation
of~$P/Q$.\end{proof}

This theorem is brought to life by the examples in the following section, which
in turn depend on the particular classical invariants given in the
appendix.  There are variations on this construction in the presence of other
canonically defined ideals in~$A$.  The singularities $\tilde E_8$ yield to
such a variation.

\section{Examples}\label{examples}
The following examples are typical. The first three were studied by
Saito~\cite{s}.
\subsection{The simple elliptic singularities $\tilde E_6$}\label{esix}
Following~\cite{csy}, we write these singularities in ${\mathbb C}^3$ in the
form
$$V_t=\{x^3+y^3+z^3+txyz=0\}\quad\mbox{for a parameter $t$ with }t^3+27\not=0$$
and present the moduli algebra as
$$\textstyle A_t={\mathbb C}\langle 1,x,y,z,yz,zx,xy,xyz\rangle
  \quad\mbox{where }
  x^2=-\frac{t}{3}yz,\; y^2=-\frac{t}{3}zx,\; z^2=-\frac{t}{3}xy.$$
Then, the maximal ideal is ${\mathfrak m}=(x,y,z)$ and we may write
$${\mathfrak m}/{\mathfrak m}^2={\mathbb C}\langle x,y,z\rangle
  \quad\mbox{and}\quad{\mathfrak m}^3={\mathbb C}\langle xyz\rangle.$$
Multiplication
$\bigodot^3{\mathfrak m}/{\mathfrak m}^2\to{\mathfrak m}^3$ has the following
effect:--
$$\textstyle x^3\mapsto -\frac{t}{3}xyz\qquad
             y^3\mapsto -\frac{t}{3}xyz\qquad
             z^3\mapsto -\frac{t}{3}xyz\qquad
             xyz\mapsto xyz$$
and all other monomials are sent to zero. As a polynomial with respect to
the dual basis $X,Y,Z$ of $({\mathfrak m}/{\mathfrak m}^2)^*$, this is
\begin{equation}\label{new}F_t(X,Y,Z)=tX^3+tY^3+tZ^3-18XYZ\end{equation}
up to scale.  Notice that this polynomial does not define the original
variety~$V_t$.  For example, it is singular if and only if $t=0$ or
$t^3=216$.  This already distinguishes these values of $t$ as special, in
agreement with Chen, Seeley, and Yau~\cite{csy} who observe that the Lie
algebra of derivations of $A_t$ jumps dimension for these values of~$t$.
The invariant theory of homogeneous cubic polynomials in three variables in
presented in the appendix.  In particular, there is an invariant $J$ of degree
4 and an invariant $K$ of degree~6.  A computation (easily carried out with
computer algebra, as explained in the appendix) gives
\begin{equation}\label{esixj}
\frac{J^3-6 K^2}{J^3}=-\frac{t^3(t^3-216)^3}{1728(t^3+27)^3}.\end{equation}
The expression on the right is the $j$-invariant (with Saito's
normalisation~\cite{s}) associated to the elliptic curve defined by the
original equation
\begin{equation}\label{original}x^3+y^3+z^3+txyz=0\end{equation}
in~${\mathbb{CP}}_2$.  (In the appendix, we shall explain how to compute the
$j$-invariant for any non-singular cubic without prior normalisation.  It is
interesting to note that the elliptic curves defined in ${\mathbb{CP}}_2$ by
(\ref{new}) and (\ref{original}) have reciprocal $j$-invariants.)
This corrects Saito's computation in \cite{s} and, multiplying by $-1728$,
gives the invariant found by Chen, Seeley, and Yau~\cite{csy}.  Once this
invariant is known it is easy to see that it is complete: the co\"ordinate
changes
\begin{equation}\label{esixfirstcoordchange}
x\mapsto\omega x\qquad y\mapsto y\qquad z\mapsto z\end{equation}
and
\begin{equation}\label{esixsecondcoordchange}
x\mapsto x+y+z\qquad y\mapsto\omega x+\omega^2y+z\qquad
y\mapsto\omega^2x+\omega y+z,\end{equation}
where $\omega^3=1$, preserve the form of $V_t$ but induce
\begin{equation}\label{esixtchange}
t\mapsto\omega t\quad\mbox{and}\quad t\mapsto\frac{3(6-t)}{3+t}\end{equation}
as regards the parameter~$t$.  These two substitutions (\ref{esixtchange})
generate the changes in $t$ that preserve~(\ref{esixj}).  In fact, Chen,
Seeley, and Yau~\cite{csy} show that, along with simply permuting and scaling
$(x,y,z)$, the substitutions (\ref{esixfirstcoordchange}) and
(\ref{esixsecondcoordchange}) generate the only co\"ordinate changes preserving
the form of~$V_t$.

\subsection{The simple elliptic singularities $\tilde E_7$}\label{eseven}
Following~\cite{sy}, we write these singularities in ${\mathbb C}^3$ in the
form
$$V_t=\{x^4+tx^2y^2+y^4+z^2=0\}\quad\mbox{for a parameter $t$ with }t^2\not=4$$
and present the moduli algebra as
$$\textstyle A_t={\mathbb C}\langle 1,x,y,x^2,xy,y^2,x^2y,xy^2,x^2y^2\rangle
  \quad\mbox{where}\quad
  x^3=-\frac{t}{2}xy^2,\; y^3=-\frac{t}{2}x^2y.$$
Then, the maximal ideal is ${\mathfrak m}=(x,y)$ and we may write
$${\mathfrak m}/{\mathfrak m}^2={\mathbb C}\langle x,y\rangle
  \quad\mbox{and}\quad{\mathfrak m}^4={\mathbb C}\langle x^2y^2\rangle.$$
Multiplication
$\bigodot^3{\mathfrak m}/{\mathfrak m}^2\to{\mathfrak m}^4$ has the following
effect:--
$$\textstyle x^4\mapsto -\frac{t}{2}x^2y^2\qquad
             x^3y\mapsto 0\qquad
             x^2y^2\mapsto x^2y^2\qquad
             xy^2\mapsto 0\qquad
             y^4\mapsto -\frac{t}{2}x^xy^2.$$
As a polynomial with respect to the dual basis $X,Y$ of
$({\mathfrak m}/{\mathfrak m}^2)^*$, this is
\begin{equation}\label{esevennew}tX^4-12X^2Y^2+tY^4,\end{equation}
up to scale.
The invariant theory of homogeneous quartic polynomials in two variables in
presented in the appendix.  In particular, there is an invariant $J$ of degree
2 and an invariant $K$ of degree~3.  A computation (easily carried out with
computer algebra, as explained in the appendix) gives
\begin{equation}\label{esevenj}
\frac{J^3}{6K^2}=\frac{(12+t^2)^3}{108(t^2-4)^2}.\end{equation}
The expression on the right is the $j$-invariant computed by Saito~\cite{s}.
The substitutions
$$t\mapsto -t\quad\mbox{and}\quad t\mapsto\frac{2(6-t)}{2+t}$$
leave (\ref{esevenj}) unchanged and generate the group of such substitutions.
These arise from the co\"ordinate changes
\begin{equation}\label{esevenchanges}
\left.\begin{array}l x\mapsto ix\\ y\mapsto y\end{array}\!\right\}
\quad\mbox{and}\quad
\left.\begin{array}l x\mapsto x+y\\ y\mapsto x-y\end{array}\!\right\},
\end{equation}
respectively.  It follows that the $j$-invariant is complete.  In fact, it is
not hard to show that, along with simply swopping and scaling $(x,y)$, the
substitutions (\ref{esevenchanges}) generate the only co\"ordinate changes
preserving the form of~$V_t$.

\subsection{The simple elliptic singularities $\tilde E_8$}\label{eeight}
Following~\cite{sy}, we write these singularities in ${\mathbb C}^3$ in the
form
$$V_t=\{x^6+tx^4y+y^3+z^2=0\}\quad\mbox{for a parameter $t$ with }
4t^3+27\not=0$$
and present the moduli algebra as
$$\textstyle A_t={\mathbb C}\langle 1,x,y,x^2,xy,x^3,x^2y,x^4,x^3y,x^4y\rangle
  \quad\mbox{where}\quad
  y^2=-\frac{t}{3}x^4,\; x^5=-\frac{2t}{3}x^3y.$$
There is the maximal ideal ${\mathfrak m}=(x,y)$ but also another canonically
defined ideal
$${\mathfrak n}=\{v\in A_t\mbox{ s.t. }v^4=0\}=(y,x^2).$$
Then
$${\mathfrak n}/{\mathfrak m}{\mathfrak n}={\mathbb C}\langle y,x^2\rangle
  \quad\mbox{and}\quad{\mathfrak n}^3={\mathbb C}\langle x^4y\rangle.$$
Multiplication
$\bigodot^3{\mathfrak n}/{\mathfrak m}{\mathfrak n}\to{\mathfrak n}^3$ has the
following effect:--
$$\textstyle y^3\mapsto -\frac{t}{3}x^4y\qquad
             y^2x^2\mapsto \frac{2t^2}{9}x^4y\qquad
             yx^4\mapsto x^4y\qquad
             x^6\mapsto -\frac{2t}{3}x^4y.$$
As polynomial with respect to the dual basis $Y,X$ of
$({\mathfrak n}/{\mathfrak m}{\mathfrak n})^*$, this is
$$tY^3-2t^2Y^2X-9YX^2+2tX^3,$$
up to scale.  There are no absolute invariants of a binary cubic under
${\mathrm{GL}}(2,{\mathbb C})$.  In fact, all classical invariants are
polynomials in a particular invariant $J$ of degree~4, proportional to the
discriminant .  However, ${\mathfrak n}/{\mathfrak m}{\mathfrak n}$ has a
canonical subspace ${\mathfrak m}^2/{\mathfrak m}{\mathfrak n}$ spanned by
$x^2$.  Therefore we may look for invariants of this binary cubic under the
subgroup of ${\mathrm{GL}}(2,{\mathbb C})$ consisting of upper triangular
matrices
\begin{equation}\label{actionofuppertriangular}
\left\lgroup\begin{array}cX\\ Y\end{array}\right\rgroup\mapsto
\left\lgroup\begin{array}{cc}a&b\\ 0&d\end{array}\right\rgroup
\left\lgroup\begin{array}cX\\ Y\end{array}\right\rgroup.\end{equation}
The invariant theory for this action is explained in the appendix.  In
particular, the coefficient of $X^3$ defines a linear invariant $K$ and there
is a simple quadratic invariant~$L$.  We may combine these to produce the
following absolute invariant:--
\begin{equation}\label{eeightj}
\frac{JK^2}{L^3}=\frac{4t^3}{4t^3+27}.\end{equation}
The expression on the right is the $j$-invariant computed by Saito~\cite{s}.
The replacement $y\mapsto\omega y$ for $\omega^3=1$ shows that there are no
further invariants.

\subsection{Another example}\label{another}
Consider the singularities in ${\mathbb C}^3$ given by
$$V_t=\{x^5+tx^3y^2+y^5+z^2=0\}\quad\mbox{for a parameter $t$ with }
108t^5+3125\not=0.$$
The moduli algebra may be generated by $x$ and $y$ subject to the relations
$$\textstyle x^4=-\frac{3t}{5}x^2y^2\quad
  \mbox{and}\quad y^4=-\frac{2t}{5}x^3y.$$
Iterating these relations quickly leads to
$$\begin{array}l\textstyle x^6=\frac{54}{625}x^3y^3,\qquad
  x^5y=-\frac{3t}{5}x^3y^3,\qquad x^4y^2=-\frac{18t^3}{125}x^3y^3,\\[4pt]
  x^2y^4=\frac{30t^2}{125}x^3y^3,\qquad xy^5=\frac{36t^4}{625}x^3y^3,\qquad
  y^6=-\frac{2t}{5}x^3y^3\end{array}$$
and so multiplication
$\bigodot^6{\mathfrak m}/{\mathfrak m}^2\to{\mathfrak m}^6$ is represented by
the polynomial
\begin{equation}\label{binsex}
\begin{array}r 27t^4X^6-1125tX^5Y-675t^3X^4Y^2+6250X^3Y^3\qquad\\[5pt]
{}+1125t^2X^2Y^4+108t^4XY^5-125tY^6\end{array}\end{equation}
up to scale.  The invariant theory of homogeneous sextic polynomials in two
variables in presented in the appendix.  In particular, there is an invariant
$J$ of degree 2 and an invariant $K$ of degree~4.  A computation (easily
carried out with computer algebra, as explained in the appendix) gives
\begin{equation}\label{anotherj}
\frac{J^2}{J^2-2K}=\frac{78125}{3(108t^5+3125)}\end{equation}
as an absolute invariant.  The replacement $y\mapsto\omega y$ for $\omega^5=1$
shows that there are no further invariants.

More generally, the family of singularities given by
\begin{equation}\label{twoparameterequation}
f(x,y,z)=x^5+sx^4y+tx^3y^2+y^5+z^2=0\end{equation}
has two absolute invariants
$$\frac{(3st^2-125)^2}{256s^5-1600s^3t-27s^2t^4+2250st^2+108t^5+
3125}$$
and
$$\frac{\left(\!\begin{array}l
163200s^6t^2+14800000s^5-2100000s^4t^3+5400s^3t^6\\[2pt]
\quad{}-92500000s^3t+7425000s^2t^4-52650st^7\\[2pt]
\qquad{}+116250000st^2+729t^{10}-4556250t^5+312500000\end{array}\!\right)^2}
{(256s^5-1600s^3t-27s^2t^4+2250st^2+108t^5+3125)^3}$$
derived in a similar way.

\appendix
\section*{Appendix: classical invariant theory}
We review Weyl's classical theory~\cite{w}, though most of what we shall need
was already well-known in the nineteenth century~\cite{e}.  Suppose
$P(X^1,\dots,X^n)$ is a homogeneous polynomial of degree~$N$.  Then we may
write
\begin{equation}\label{converttotensor}p(X^1,\dots,X^n)=\sum_{i,j,\dots,k=1}^n
a_{\mbox{\scriptsize$\underbrace{ij\cdots k}_{\makebox[0pt]{$N$ indices}}$}}
X^iX^j\cdots X^k\;=\;a_{ij\cdots k}X^iX^j\cdots X^k\end{equation}
for a tensor $a_{ij\cdots k}$, symmetric in its indices. In the last
expression, the summation is implicit---the Einstein summation convention
demands a summation over repeated indices.  We shall use this convention form
now on.  Fix a totally skew tensor
$\epsilon^{i\cdots j}$ with $n$ indices.  More invariantly, we may view the
tensor $a$ as an element of $\bigodot^N{\mathbb V}$ for an $n$-dimensional
vector space ${\mathbb V}$ and then $\epsilon$ is chosen
in~$\Lambda^n{\mathbb V}^*$.  Since this space is one-dimensional, $\epsilon$
is unique up to scale. In practise, we can take $\epsilon^{12\cdots n}=1$ to
fix the scale.

As an example, consider the polynomial~(\ref{original}). It is homogeneous of
degree 3 in 3 variables, classically a `ternary cubic'---see~\cite{e}.
Converting to tensor notation in accordance with~(\ref{converttotensor}), we
see that
$$a_{111}=1,\;a_{222}=1,\;a_{333}=1,\;
  a_{123}=a_{231}=a_{312}=a_{213}=a_{132}=a_{213}=t/6$$
and all other $a_{ijk}$ are zero.  Now consider
\begin{equation}\label{Jforternarycubic}J=a_{ijk}a_{lmn}a_{pqr}a_{stu}
  \epsilon^{ilp}\epsilon^{jms}\epsilon^{kqt}\epsilon^{nru}.\end{equation}
Viewing $a$ and $\epsilon$ invariantly, as elements of $\bigodot^3{\mathbb V}$
and $\Lambda^3{\mathbb V}^*$ respectively, it is clear that $J$ is independent
of any choice of basis: a summation over repeated indices is simply the
invariant contraction ${\mathbb V}\otimes{\mathbb V}^*\to{\mathbb C}$.  But
$\epsilon^{ijk}$ was chosen arbitrarily and scales by $\det(M)$ under the
action of $M\in{\mathrm {GL}}({\mathbb V}^*)$.  It follows that if we change
co\"ordinates $X^i\mapsto m^i{}_jX^j$ in the polynomial $p(X^1,\dots,X^n)$,
and recompute $J$, then the answer will be multiplied by $\det(m^i{}_j)^4$.
Using the original co\"ordinates and choosing $\epsilon^{123}=1$, we find that
$J=t(t^3-216)/54$.  Such calculations may be performed on a computer in
just a few seconds once a suitable program is written.  All the calculations in
this article were done this way using the {\sc Maple} computer algebra package.
The programs are available electronically\footnote{
ftp://ftp.maths.adelaide.edu.au/pure/meastwood/maple/README}.

Expressions such as (\ref{Jforternarycubic}), in which various
tensors are juxtaposed with all indices
paired up, are called `complete contractions'.  Weyl's first fundamental
theorem of invariant theory~\cite{w}
states that all ${\mathrm{SL}}(n,{\mathbb C)}$
invariants of any tensor $a_{ij\cdots k}$ arise as linear combinations of
complete contractions of $a_{ij\cdots k}$ and~$\epsilon^{i\cdots j}$.  More
complicated complete contractions are best built up from `partial
contractions' (or classical `covariants'---see~\cite{e}) such as
$$b_{ij}{}^{kl}=a_{pqi}a_{rsj}\epsilon^{prk}\epsilon^{qsl}\qquad
   c_{ijk}=b_{ij}{}^{pq}a_{pqk}\qquad
d^{ijk}=b_{pq}{}^{ir}a_{rst}\epsilon^{ptk}\epsilon^{qjs}.$$
For example,
$$J=b_{ij}{}^{kl}b_{kl}{}^{ij}\quad\mbox{and}\quad
K=c_{ijk}d^{ijk}$$
reproduces $J$ as in (\ref{Jforternarycubic}) and gives a new invariant $K$
homogeneous of degree~6.  It is a classical theorem~\cite{e} that these $J$ and
$K$ freely generate the ring of all ${\mathrm{SL}}(3,{\mathbb C})$-invariants
of a ternary cubic.
\begin{theorem}\label{jtheorem} The ratio
$$j=\frac{J^3}{J^3-6K^2}$$
associated to any homogeneous polynomial p(x,y,z) of degree 3 is independent
of choice of co\"ordinates.  When the denominator is non-zero, this is the
$j$-invariant of the elliptic curve $\{p(x,y,z)=0\}$ in~${\mathbb{CP}}_2$.
\end{theorem}
\begin{proof}Since numerator and denominator are homogeneous of the same
degree, namely~12, a change of co\"ordinates effected by a non-singular matrix
$M$ scales both by $\det(M)^{12}$.  Though evident by construction, it is also
easily checked by direct computation with a computer.  In any case, $j$ is
co\"ordinate-independent.  To show that it is the $j$-invariant, it suffices
to verify this for one of the canonical forms of a non-singular cubic.  For
example, a computer calculation shows
$$f(x,y,z)=zy^2-x(x-z)(x-\lambda z)\Longrightarrow
j=\frac{4}{27}\frac{(\lambda^2-\lambda+1)^3}{\lambda^2(\lambda-1)^2}.$$
This is the well-known formula (normalised as in~\cite{s}).
\end{proof}
Several canonical forms (including Weierstra\ss) are included in the
computer program `ternary\_cubic' but any polynomial can be treated: for
example,
$$p(x,y,z)=x^3+x^2y-4z^3+xyz-xz^2+xy^2\Longrightarrow
j=\frac{357911}{120545280}.$$
This would be difficult to compute if one first had to transform it into a
canonical form but using Theorem~\ref{jtheorem} and a computer takes just a few
seconds.  Similarly,
$$p(x,y,z)=x^3+y^3+z^3+txyz\Longrightarrow
j=-\frac{t^3(t^3-216)^3}{1728(t^3+27)^3}.$$
This corrects the calculation of $j$ in~\cite{s}.  On the other hand, it is
the reciprocal that computes the $j$-invariant via the moduli algebra
and~(\ref{new}):--
$$p(X,Y,Z)=tX^3+tY^3+tZ^3-18XYZ\Longrightarrow
\frac{J^3-6K^2}{J^3}=-\frac{t^3(t^3-216)^3}{1728(t^3+27)^3}.$$
This completes the
invariant theory of the ternary cubic used in \S\ref{esix}.

In~\S\ref{eseven}, it was the invariant theory of a homogeneous quartic in two
variables that was used. The classical theory of `binary quartics' is given in
Elliott~\cite{e}.  The ring of invariants is freely generated by
$$J=b_{ij}{}^{kl}b_{kl}{}^{ij}\quad\mbox{and}\quad
K=b_{ij}{}^{kl}b_{kl}{}^{mn}b_{mn}{}^{ij},\quad\mbox{where }
b_{ij}{}^{kl}=a_{ijpq}\epsilon^{pk}\epsilon^{ql}.$$
These are computed in the program `binary\_quartic'. In particular,
$$p(x,y)=x(x-y)(x-\lambda y)y\Longrightarrow \frac{4J^3}{4J^3-24K^2}=
\frac{4}{27}\frac{(\lambda^2-\lambda+1)^3}{\lambda^2(\lambda-1)^2}$$
so this ratio is the $j$-invariant in general. Then,
$$p(x,y)=x^4+tx^2y^2+y^4\Longrightarrow j=\frac{(12+t^2)^3}{108(t^2-4)^2},$$
verifying Saito's result~\cite{s}.
On the other hand, via the moduli algebra and (\ref{esevennew}):--
$$p(X,Y)=tX^4-12X^2Y^2+tY^4\Longrightarrow
\frac{J^3}{6K^2}=\frac{(12+t^2)^3}{108(t^2-4)^2},$$
as claimed in~\S\ref{eseven}.

The invariant theory of the binary sextic was used in~\S\ref{another}.  As
discussed in~\cite{e}, the ring of invariants has five generators $J$, $K$,
$L$, $M$, and $N$ of degrees 2, 4, 6, 10, and 15, respectively.  To define
them, first consider the covariants
$$\begin{array}l b^{ijklmn}=a_{pqrstu}
\epsilon^{pi}\epsilon^{qj}\epsilon^{rk}\epsilon^{sl}\epsilon^{tm}\epsilon^{un}
\qquad c_{ij}{}^{kl}=a_{ijpqrs}b^{klpqrs}\\[5pt]
d_{ij}=a_{ijpqrs}c_{tu}{}^{rs}\epsilon^{pt}\epsilon^{qu}
\qquad f_{ij}=c_{ij}{}^{pq}d_{pq}
\qquad g_{ij}=c_{ij}{}^{pq}f_{pq}.\end{array}$$
Then,
$$\begin{array}l J=c_{ij}{}^{ij}\qquad
K=c_{ij}{}^{kl}c_{kl}{}^{ij}\qquad
L=c_{ij}{}^{kl}c_{kl}{}^{mn}c_{mn}{}^{ij}\\[5pt]
M=b^{pqrstu}d_{pq}d_{rs}d_{tu}\qquad
N=\epsilon^{qr}\epsilon^{st}\epsilon^{up}d_{pq}f_{rs}g_{tu}.\end{array}$$
These are computed in the program `binary\_sextic'. There is just one relation
$$N^2+\frac{1}{1458}J^{15}
             -\frac{7}{486}J^{13}K
             +\frac{13}{108}J^{11}K^2
             +\cdots
             +\frac{1}{8}KLM^2
             +\frac{1}{16}M^3=0$$
given explicitly and verified in binary\_sextic. For the purposes
of~\S\ref{another}, the program checks (\ref{anotherj}) for the
polynomial~(\ref{binsex}).

In the case of (\ref{twoparameterequation}), there is an identity
$J^3+3JK-10L=0$ and so we are obliged to use $M$ if we want to construct two
absolute invariants in accordance with Theorem~\ref{construction}. The two
invariants listed for this case are
$$\frac{3}{5}\frac{J^2}{J^2-2K}\quad\mbox{and}
\quad 759375\frac{M^2}{(J^2-2K)^5}.$$

Finally we discuss the invariant theory behind~\S\ref{eeight}.
The basic invariant of a binary cubic $a_{ijk}X^iX^jX^k$ is
$$J=c_i{}^jc_j{}^i\quad\mbox{where }
c_i{}^j=a_{ipq}\epsilon^{jk}\epsilon^{pl}\epsilon^{qm}a_{klm}.$$
It is $-2/27$ times the discriminant. The ring of
${\mathrm{SL}}(2,{\mathbb C})$-invariants consists of polynomials in $J$ so
there is no possibility of applying Theorem~\ref{construction} directly.
Instead, as suggested in~\S\ref{eeight}, we may look for invariants under the
action (\ref{actionofuppertriangular}) of the upper triangular matrices.
Consider the covector $e^i$ with $e^1=1$ and $e^2=0$. It is distinguished by
being preserved up to scale by~(\ref{actionofuppertriangular}). Therefore, it
is another ingredient that may be used in complete contractions such as
$$K=a_{ijk}e^ie^je^k\qquad\mbox{and}\qquad
L=a_{ijk}a_{lmn}\epsilon^{il}\epsilon^{jm}e^ke^n$$
to produce polynomials that transform by a character
under~(\ref{actionofuppertriangular}). Indeed, it is shown in~\cite{f} (as a
very special case of an invariant theory for tensor representations of
parabolic subgroups of the classical groups) that all such polynomials arise as
linear combinations of complete contractions like this. In any case, the
quotient $JK^2/L^3$ is not only of polynomials of equal homogeneity but also of
the same degree in $\epsilon$, namely~$6$. It is, therefore, an absolute
invariant of a cubic under~(\ref{actionofuppertriangular}), as required.

\end{document}